\documentstyle[12pt]{article}

\begin{document}

\def\myone{\hbox{\bf 1}}
\def\myzero{\hbox{\bf 0}}
\def\arctg{\operatorname{arctg}}
\def\supp{\operatorname{supp}}
\def\cond{\operatorname{cond}}
\def\Tr{\operatorname{Tr}}
\def\tg{\operatorname{tg}}
\def\th{\operatorname{th}}
\def\pd#1#2{\frac{\partial#1}{\partial#2}}
\def\sumpr{\mathop{{\sum}'}}

\title{Error Autocorrection in Rational Approximation and Interval
Estimates} 
\author{Grigori L. Litvinov 
\thanks{Partly supported by the Fields Institute for Research
in Mathematical Sciences (Toronto, Canada).}}      

\date{}

\maketitle

\begin{center}
{\it Independent University of Moscow

B. Vlasievskii per., 11, Moscow 121002, Russia

E-mail: islc@dol.ru and glitvinov@mail.ru}
\end{center}

\bigskip
\begin{abstract}
The error autocorrection effect means that in a calculation all the
intermediate errors compensate each other, so the final result
is much more accurate than the intermediate results. In this case
standard interval estimates (in the framework of interval analysis 
including the so-called a posteriori interval analysis of Yu.
Matijasevich) are too pessimistic. We shall discuss a very strong
form of the effect which appears in rational approximations to
functions.

The error autocorrection effect occurs in all efficient methods of
rational approximation (e.g., best approximations, Pad\'e
approximations, multipoint Pad\'e approximations, linear and nonlinear
Pad\'e-Chebyshev approximations, etc.), where very significant errors
in the approximant coefficients do not affect the accuracy of this
approximant. The reason is that the errors in the
coefficients of the rational approximant are not distributed in an
arbitrary way, but form a collection of coefficients for a new
rational approximant to the same approximated function. The
understanding of this mechanism allows to decrease the approximation
error by varying the approximation procedure depending on the form of
the approximant. Results of computer experiments are presented.

The effect of error autocorrection indicates that variations of an
approximated function under some deformations of rather a general
type may have little effect on the corresponding rational approximant
viewed as a function (whereas the coefficients of the approximant can
have very significant changes). Accordingly, while deforming a
function for which good rational approximation is possible, the
corresponding approximant's error can rapidly increase, so the
property of having good rational approximation is not stable under
small deformations of the approximated functions. This property is
``individual'', in the sense that it holds for specific functions.



\end{abstract}

\section{Introduction}
The paper contains a brief survey description of results published in 
\cite{1}--\cite{7} (but not all of them) and a new material on a 
general form of the error autocorrection effect and its relation to 
interval analysis in the spirit of \cite{8}--\cite{10}. Some old 
results are presented in a new form. The treatment is partly 
heuristic and it is based on computer experiments.

The error autocorrection effect means that in a calculation all the
intermediate calculating errors compensate each other, so the final result
is much more accurate than the intermediate results. In this case
standard interval estimates (e.g., in the framework of Yu. 
Matijasevich's a posteriori interval analysis) 
are not realistic, they are too pessimistic. 
The error autocorrection effect appears in some popular numerical
methods, e.g., in the least squares method. In principle this
effect is not new for experts in interval computations.
We shall discuss a very strong form of the effect.

For the sake of simplicity let us suppose that we calculate values of
a real smooth function $z=F(y_1, \cdots, y_n)$ of real variables 
$y_1, \cdots, y_n$, and suppose that all round-off errors are
negligible with respect to input data errors. In this case the error 
$\Delta F$ of the function $F$ can be evaluated by the formula
$$
\Delta F=\sum_{i=1}^N \pd F{y_i}\cdot \Delta y_i+r,
$$
where $r$ is negligible. The sum $\sum_{i=1}^N \pd F{y_i}\cdot \Delta y_i$
in this formula can be treated as a scalar product of the gradient
vector $\{\pd F{y_i}\}$ and the vector of errors
$\{\Delta y_i\}$. For all the standard interval methods the best
possible estimate for $\Delta F$ is given by the formula
$$
\vert \Delta F\vert \leq \sum_{i=1}^N \vert\pd F{y_i}\vert
\cdot \vert\Delta y_i\vert.
$$
This estimate is too pessimistic if both the vectors
$\{\pd F{y_i}\}$ and $\{ \Delta y_i\}$ are not small but 
the scalar product is small enough, so these vectors are almost
orthogonal. That is the case of the error autocorrection effect. 
To calculate values of a function we usually use four arithmetic
operations (applying to arguments and constants) and this leads
to a rational approximation to the calculated function.

The author came across the phenomenon of error autocorrection at the end of
seventies while developing nonstandard algorithms for computing elementary
functions on small computers. It was required to construct rational
approximants of the form
\begin{equation}
R(x)={a_0+a_1x+a_2x^2+\dots +a_nx^n\over{b_0+b_1x+b_2x^2+\dots +b_mx^m}}
\end{equation}
to certain functions of one variable $x$ defined on finite segments of the
real line. For this purpose a simple method (described in \cite{1} and
below) was used: the method allows to determine the family of
coefficients $a_i$, $b_j$ of the approximant~(1) as a solution of a
certain system of linear algebraic equations. These systems
turned out to be ill conditioned, i.e., the problem of
determining the coefficients of the approximant is, generally speaking,
ill-posed and small perturbations of the approximated function $f(x)$ or 
calculation errors lead to very significant errors 
in the values of coefficients. Nevertheless, the method ensures a
paradoxically high quality of the obtained approximants whose errors are
close to the best possible \cite{1, 2}.

For example, for the function $\cos x$ the approximant of the form~(1) on
the segment $[-\pi /4,\pi /4]$ obtained by the method mentioned above for
$m=4$, $n=6$ has the relative error equal to
$0.55\cdot 10^{-13}$, and the best
possible relative error is $0.46\cdot 10^{-13}$ \cite{11}. The corresponding
system of linear algebraic equations has the condition number of order
$10^9$. Thus we risk to lose 9 accurate decimal digits in the
solution  because of calculation errors. Computer experiments
show that this is a serious risk. The method mentioned above
was implemented as a Fortran program. The calculations were
carried out with double precision (16 decimal positions) by means of
two different computers. These computers are very similar in their
architecture, but when passing from one computer to another the system of
linear equations and the computational process are perturbed because of
calculation errors, including round-off errors. As a result, the
coefficients of the approximant mentioned above to the function $\cos x$
experience a perturbation already at the
sixth--ninth decimal digits. But the error in the rational approximant
itself remains
invariant and is $0.4\cdot 10^{-13}$ for the absolute error and $0.55\cdot
10^{-13}$ for the relative error. The same thing happens for approximants
of the
form~(1) to the function $\arctan x$ on
the segment $[-1,1]$ obtained by the method mentioned above
for $m=8$, $n=9$ the relative error is $0.5\cdot 10^{-11}$
and does not change while passing from one computer to another although the
corresponding system of linear equations has the condition number of
order $10^{11}$, and the coefficients of the approximant experience a
perturbation with a relative error of order~$10^{-4}$.

Thus the errors in the numerator and the denominator of a rational 
approximant
compensate each other. The effect of error autocorrection is connected with
the fact that the errors in the coefficients of a rational approximant are 
not
distributed in an arbitrary way, but form the coefficients of a new
approximant to the approximated function. It can be easily understood that 
all
standard methods of interval arithmetic (see, for example, \cite{6}, 
 \cite{8}--\cite{10})
do not allow us to take into account this effect
and, as a result, to estimate the error in the rational approximant
accurately (see section 12 below).

Note that the application of standard procedures known
in the theory of ill-posed problems results in this case in losses
in accuracy. For example, if one applies the regularization method,
then two thirds of the accurate figures are lost \cite{12};
in addition, the amount of calculations increases rapidly. The
matter of import is that the exact solution of the system of equations
in the present case is not the ultimate goal; the aim is to construct an
appro\-ximant which is precise enough. This approach allows to
``rehabilitate'' (i.e., to justify) and to simplify a number of
algorithms intended for the construction of the approximant, to obtain 
(without additional transforms) approximants in a form which is convenient 
for applications.

Professor Yudell L.~Luke kindly drew the author's attention to his
papers \cite{13, 14} where the effect of error autocorrection for the classical
Pad\' e approximants was revealed and was explained at a heuristic level.
The method mentioned above leads to the linear
Pad\' e--Chebyshev approximants if the calculation errors are
ignored.

In the present paper the error autocorrection mechanism is considered 
for quite a general situation (including linear methods for the 
construction of rational approximants and nonlinear generalized Pad\'e 
approximations). The efficiency of the construction algorithms used for 
rational approximants is due to the error autocorrection effect (at least 
in the case when the number of coefficients is large enough).

Our new understanding of the error autocorrection mechanism allows us, 
to some extent, to control calculation errors by changing
the construction procedure depending on the form of the approximant.
It is shown that the appearance of a control parameter allowing to take 
into account the error autocorrection mechanism ensures the decrease of 
the calculation errors in some cases.

Construction methods for linear and nonlinear Pad\'e--Chebyshev 
approximants involving the computer algebra system REDUCE (see \cite{4}) 
are also briefly described. Computation results characterizing the 
comparative precision of these methods are given. With regard to the 
error autocorrection phenomenon the effect described in \cite{2} and 
connected with the fact that a small variation of an approximated 
function can lead to a sharp decrease in accuracy of the Pad\'e--Chebyshev 
approximants is analyzed. 

The error autocorrection effect occurs not only in rational 
approximation but appears in some other cases (approximate solutions 
of linear differential equations, the method of least squares etc.). In this 
more general situation relations between the error autocorrection effect and 
standard methods of interval analysis are also discussed in this paper.

The author is grateful to Y.L.~Luke, B.S. Dobronets, and S.P Shary for 
stimulating discussions and comments. The author wishes to express his 
thanks to I.~A.~Andreeva, A.~Ya.~Rodionov and V.~N.~Fridman who 
participated in the programming and organization of computer experiments.

\section{Error autocorrection in rational approximation}

Let  $ \{\varphi_0,\varphi_1,\dots ,\varphi_n\}$  and
$\{\psi_0,\psi_1,\dots ,\psi_m\}$ be two collections consisting of linearly
independent functions of the argument $x$ belonging to some (possibly
multidimensional) set $X$. Consider the problem of constructing
an approximant of the form
\begin{equation}
R(x)={a_0\varphi_0+a_1\varphi_1+\dots +a_n\varphi_n\over
b_0\psi_0+b_1\psi_1+\dots +b_m\psi_m}
\end{equation}
to a given function $f(x)$ defined on $X$. If $X$ coincides with a real line
segment $[A,B]$ and if $\varphi_k=x^k$ and $\psi_k=x^k$ for
all $k$, then the expression~(2) turns out to be a rational function
of the form~(1)
(see the Introduction). It is clear that expression~(2) also gives a
rational function in the case when we take Chebyshev polynomials
$T_k$ or, for example, Legendre, Laguerre, Hermite, etc.
polynomials as $\varphi_k$ and $\psi_k$.

Fix an abstract construction method for an approximant of the form (2)
and consider the problem of computing the coefficients $a_i$, $b_j$.
Quite often this problem is ill-conditioned (ill-posed). For example, 
the problem of computing coefficients for best rational approximants 
(including polynomial approximants) for high degrees of the numerator 
or the denominator is ill-conditioned.

The instability with respect to the calculation error can be related
both to the abstract construction method of approximation (i.e.,
with the formulation of the problem) and to the particular algorithm
implementing the method. The fact that the problem of computing
coefficients for the best approximant is ill-conditioned is related
to the formulation of this problem. This is also valid for other 
construction methods for rational approximants with a sufficiently 
large number of coefficients. But an unfortunate choice of the
algorithm implementing a certain method can aggravate troubles
connected with ill-conditioning.

 Let the coefficients $a_i$, $b_j$ give an exact or an approximate 
solution of this problem, and let the $\tilde a_i$, $\tilde b_j$ give 
another approximate solution obtained in the same way. Denote by 
$\Delta a_i$, $\Delta b_j$ the absolute errors of the coefficients, 
i.e., $\Delta a_i=\tilde a_i-a_i$, $\Delta b_j=\tilde b_j-b_j$; 
these errors arise due to perturbations of the approximated function 
$f(x)$ or due to calculation errors. Set
$$
\begin{array}{rclrcl}
P(x)&=&\sum_{i=0}^na_i\varphi_i,&\quad Q(x)&=&\sum_{j=0}^mb_j\psi_j,\\
\Delta P(x)&=&\sum_{i=0}^n\Delta a_i\varphi_i, &\quad \Delta Q(x)&=&
\sum_{j=0}^m\Delta b_j\psi_j,\\
\widetilde{P}(x)&=&P+\Delta P,&\quad \widetilde{Q}(x)&=&Q+\Delta Q.
\end{array}
$$
It is easy to verify that the following exact equality is valid:
\begin{equation}
{P+\Delta P\over Q+\Delta Q}-{P\over Q}={\Delta Q\over Q}\Bigg({\Delta
P\over
\Delta Q}-{P\over Q}\Bigg). 
\end{equation}
 As mentioned in the Introduction, the fact that the problem of calculating
coefficients is ill-conditioned can nevertheless be accompanied by
high accuracy of the
approximants obtained. This means that the approximants $P/Q$ and
$\widetilde{P}/\widetilde{Q}$ are close to the approximated function and,
therefore, are close to each other, although the coefficients of these
approximants differ greatly. In this
case the relative error $\Delta Q/ \widetilde{Q}=\Delta Q/(Q+\Delta Q)$ of
the denominator considerably exceeds in absolute value the left-hand
side of equality~(3). This is possible only in the case when the
difference $\Delta P/\Delta Q-P/Q$ is small, i.e., the function
$\Delta P/\Delta Q$ is close to $P/Q$, and, hence, to the approximated
function. Thus the function $\Delta P/\Delta Q$ will be called the
{\it error approximant}.
For a special case, this concept was actually introduced in \cite{13}.
For "efficient" methods, the error approximant
provides indeed a good
approximation for the approximated function and, thus, $P/Q$ and
$\widetilde{P}/\widetilde{Q}$ differ from each other by a product of
small quantities in the right-hand side of~(3). 

Usually the following {\bf ``uncertainty relation''} is valid:
$$
\Delta\Bigg({P\over Q} \Bigg) =\delta Q \cdot\Bigg({\Delta P\over\Delta Q}
-{P\over Q}\Bigg)\approx\delta Q \cdot \Bigg({\Delta P\over\Delta Q}
-f\Bigg)\sim\varepsilon,
$$
where $\delta Q ={\Delta Q\over Q+\Delta Q}$ is the relative error of
 the denominator $Q$, the difference ${\Delta P\over\Delta Q}-f$ is the 
absolute error of the error approximate ${\Delta P\over\Delta Q}$
to the function $f$, and $\varepsilon$ is the absolute ``theoretical'' 
error of our method; the argument $x$  can be treated as fixed.

Usually the approximation function $f(x)$  is treated  as an element of 
a Banach space with a norm $\| \cdot \|$ and the absolute error $\Delta$
of the approximant ~(2) is defined as $\Delta=\| f-{P\over Q}\|$; so its 
relative error $\delta$ is defined as $\delta=\|(f-{P\over Q})/f\|$ or 
$\delta=\|(f-{P\over Q})/{P\over Q}\|$.
In what follows, we shall consider the space $C[A,B]$ of all continuous 
functions defined on the real line segment $[A,B]$ with the norm
$$\| f(x)\|= \max_{A\le x\le B}\vert f(x)\vert.$$
Below we discuss examples of the error autocorrection effect for linear 
and nonlinear methods of rational approximation.

\section{Error autocorrection for linear methods in  rational approximation}

Several construction methods for approximants of the form~(2)
are connected with solving systems of linear algebraic equations. This
procedure can lead to a large error if the corresponding matrix is
ill-conditioned. Consider an arbitrary system of linear algebraic equations
\begin{equation}
Ay=h,
\end{equation}
where $A$ is a given square matrix of order $N$ with components
$a_{ij}$ $(i,j=1,\dots ,N)$, $h$ is a given vector column with components
$h_i$, and $y$ is an unknown vector column with
components $y_i$. Define the vector norm by the equality
$$
\Vert y\Vert =\sum_{i=1}^N\vert x_i\vert
$$
(this norm is more convenient for calculations than
$\sqrt{x_1^2+\dots +x_N^2}$). Then the matrix norm is determined by the
equality
$$
\| A\| =\max_{\| y\| =1}\| Ay\| =\max_{1\le j\le N}
\sum_{i=1}^N\| a_{ij}\|.
$$
If a matrix $A$ is nonsingular, then the quantity
\begin{equation}
cond (A)=\|A\|\cdot\|A^{-1}\| 
\end{equation}
is called the {\it condition number} of the matrix $A$ (see, for example,
\cite{15}). Since $y=A^{-1}h$, we see that the absolute error $\Delta y$
of the vector $y$ is connected with the absolute error of the vector $h$
by the relation $\Delta y=A^{-1}\Delta h$, whence
$$
\Vert\Delta y\Vert\leq\Vert A^{-1}\Vert\cdot\Vert\Delta h\Vert
$$
and
$$ 
\Vert\Delta y\Vert /\Vert y\Vert\leq\Vert A^{-1}\Vert\cdot
(\Vert h\Vert /\Vert y\Vert)(\Vert\Delta h\Vert /\Vert h\Vert).
$$
Taking into account the fact that
$\Vert h\Vert\leq\Vert A\Vert\cdot\Vert y\Vert$, we finally obtain
\begin{equation}
\Vert\Delta y\Vert /\Vert y\Vert\leq\Vert A\Vert\cdot\Vert A^{-1}
\Vert\cdot\Vert\Delta h\Vert /\Vert h\Vert ,
\end{equation}
i.e., the relative error of the solution $y$ is estimated via the relative
error of the vector $h$ by means of the condition number. It is clear
that~(6) can turn into an equality. Thus, if the condition number
is of order $10^k$, then, because of round--off errors in $h$, we
can lose $k$ decimal digits of $y$.

The contribution of the error of the matrix $A$ is evaluated similarly. 
Finally, the dependence of $cond (A)$ on the choice of a norm is weak.
A method of rapid estimation of the condition number is described
in \cite{15}, Section 3.2. 

Let an abstract construction method for the approximant of the
form (2)
be linear in the sense that the coefficients of the approximant can be
determined from a homogeneous system of linear algebraic equations. The
homogeneity condition is connected with the fact that, when multiplying
the numerator and the denominator of fraction~(2) by the same nonzero
number, the approximant~(2) does not change. Denote by $y$ the vector whose
components are the coefficients $a_0, a_1,\dots ,a_n$,
$b_0,b_1,\dots ,b_m$.
Assume that the coefficients can be obtained from the homogeneous system
of equations
\begin{equation}
Hy=0,
\end{equation}
where $H$ is a matrix of dimension $(m+n+2)\times (m+n+1)$.

The vector $\tilde y$ is an approximate solution of system~(7) if the
quantity $\Vert H\tilde y\Vert$ is small. If $y$ and $\tilde y$ are
approximate solutions of system~(7), then the vector $\Delta y=\tilde y-y$
is also an approximate solution of this system since
$\Vert H\Delta y\Vert=\Vert H\tilde y-Hy\Vert\leq\Vert H\tilde y\Vert
+\Vert Hy\Vert$. Thus it is natural to assume that the function
$\Delta P/\Delta Q$ corresponding to the solution $\Delta y$ is
an approximant to $f(x)$. It is clear that the order of the residual
of the approximate solution $\Delta y$ of system~(7), i.e., of the
quantity $\Vert H\Delta y\Vert$, coincides with the order of the largest
of the residuals of the approximate solutions $y$ and $\tilde y$. For a fixed
order of the residual the increase of
the error $\Delta y$ is compensated by the fact that $\Delta y$
satisfies the system of equations~(7)
with greater ``relative'' accuracy, and
the latter, generally speaking, leads to the increase of
accuracy of the error approximant.

To obtain a certain solution of system~(7), one usually adds to this system
a normalization condition of the form
\begin{equation}
\sum_{i=0}^n\lambda_ia_i+\sum_{j=0}^m\mu_jb_j=1,
\end{equation}
where $\lambda_i$, $\mu_j$ are numerical coefficients. As a rule, the
relation $b_0=1$ is taken as the normalization condition (but this
is not always successful with respect to minimizing the calculation
errors).

Adding equation~(8) to system~(7), we obtain a nonhomogeneous system of
$m+n+2$ linear algebraic equations of type~(4). If the approximate
solutions $y$ and $\tilde y$ of system~(7) satisfy condition~(8), then
the vector $\Delta y$ satisfies the condition
\begin{equation}
\sum_{i=0}^n\lambda_i\Delta a_i+\sum_{j=0}^m\mu_j\Delta b_j=0.
\end{equation}

Of course, the above reasoning is not very rigorous; for each
specific construction method for approximations it is necessary to carry
out some additional analysis. More accurate arguments are given below 
for the linear and nonlinear Pad\'e--Chebyshev approximants. The presence 
of the error autocorrection mechanism described above is also verified 
by numerical experiments (see below).

It is clear, that classical Pad\'e approximations, multipoint Pad\'e 
approximations, linear  generalized Pad\'e approximations in sense of 
\cite{16} (e.g., linear Pad\'e--Chebyshev approximations) give us good 
examples of linear methods in rational approximation. From our point 
of view, methods of the best approximations can be treated as linear. 
The thing is that  the coefficients of the best Chebyshev approximant 
satisfy a system of linear algebraic equations and are computed as 
approximate solutions of this system on the last step of the iteration 
process in algorithms of Remez's type (see \cite{7, 17} for details).
Thus, the construction methods for the best rational approximants can be
regarded as linear. At least for some functions (say, for
$\cos((\pi /4)x)$, $-1\leq x\leq 1$) the linear and the nonlinear
Pad\'e--Chebyshev approximants are very close to the best ones in the 
sense of the relative and the absolute errors, respectively. The results 
that arise when applying calculation algorithms for Pad\'e--Chebyshev 
approximants can be regarded as approximate solutions of the system 
which determines the best approximants.
Thus the presence of the effect of error autocorrection for Pad\'e--Chebyshev
approximants gives an additional argument in favor of the conjecture
that this effect also takes place for the best approximants.

Finally, note that the basic relation~(3) becomes meaningless if one seeks
an approximant in the form
$a_0\varphi_0+a_1\varphi_1+\dots +a_n\varphi_n$, i.e., the denominator
in~(2) is reduced to $1$. However, in this case the effect of error
autocorrection (although much weakened) is also possible; this is
connected with
the fact that the errors $\Delta a_i$ approximately satisfy certain
relations. Such a situation can arise when using the least
squares method.

\section{Linear Pad\'e--Chebyshev \quad approximations and the PADE program}

Let us start to discuss a series of examples. Consider the approximant 
of the form ~(1)
\begin{equation}
R_{m,n}(x)={a_0+a_1x+a_2x^2+\dots +a_nx^n\over{b_0+b_1x+b_2x^2+\dots 
+b_mx^m}} 
\end{equation}
to a function $f(x)$ defined on the segment $[-1,1]$. The absolute 
error function 
$$
\Delta (x)=f(x)-R_{m,n}(x)
$$ 
obviously has the following form:
$$
\Delta (x)=\Phi (x)/Q_m(x),
$$
where
\begin{equation}
\Phi (x)=f(x)Q_m(x)-P_n(x). 
\end{equation}
The function $R_{m,n}(x)=P_n(x)/Q_m(x)$ is called the {\it linear
Pad\'e--Chebyshev approximant to the function} $f(x)$ if
\begin{equation}
\int\limits_{-1}^1\Phi (x)T_k(x)w(x)\,dx=0,\qquad k=0,1,\dots ,m+n, 
\end{equation}
where $T_k(x)=\cos(n \arccos x) $ are the Chebyshev polynomials,
$w(x)=1/\sqrt{1-x^2}$. This concept allows a generalization to the 
case of other orthogonal polynomials (see, e.g., \cite{16}, 
\cite{18}--\cite{20}).  Approximants of this kind always
exist \cite{18}.  The system of
equations~(12) is equivalent to the following system of linear algebraic
equations with respect to the coefficients $a_i$, $b_j$:
\begin{equation}
\sum_{j=0}^mb_j\int\limits_{-1}^1{x^jT_k(x)f(x)\over\sqrt{1-x^2}}\,dx
-\sum_{i=0}^na_i\int\limits_{-1}^1{x^iT_k(x)\over\sqrt{1-x^2}}\,dx=0. 
\end{equation}

The homogeneous system~(12) can be transformed into a nonhomogeneous one by
adding a normalization condition; in particular, any of the following
relations can be taken as this condition:
\begin{eqnarray}
b_0&=&1,\\
b_m&=&1,\\
a_m&=&1.
\end{eqnarray}
In \cite{1, 2} the program PADE (in Fortran, with double precision) which
allows us to construct rational approximants by solving the system of
equations of type~(13) is briefly described.
The complete text of a certain version of this program and its detailed
description can be found in the Collection of algorithms
and programs of the Research Computer Center of the Russian Academy of
Sciences \cite{3}. For even functions the approximant is looked for in the form
\begin{equation}
R(x)={a_0+a_1x^2+\dots +a_n(x^2)^n\over b_0+b_1x^2+\dots
+b_m(x^2)^m},
\end{equation}
and for odd functions it is looked for in the form
\begin{equation}
R(x)=x{a_0+a_1x^2+\dots +a_n(x^2)^n\over b_0+b_1x^2+\dots
+b_m(x^2)^m},
\end{equation}
respectively.
The program computes the values of coefficients of the approximant, the
absolute and the relative errors $\Delta=\max_{A\leq x\leq B}
\vert\Delta (x)\vert$ and $\delta=\max_{A\leq x\leq B}\vert\Delta (x)
/f(x)\vert$, and gives the information which allows us to estimate the 
quality of the approximation (see \cite{7} and \cite{3} for details). 
Using a
subroutine, the user introduces the function defined by means of any
algorithm on an arbitrary segment~$[A,B]$, introduces the boundary points 
of this segment, the numbers $m$ and $n$, and the number of control 
parameters. In particular, one can choose the normalization condition 
of type~(14)--(16), look for an approximant in the form~(17) or~(18) and 
so on. The change of the variable reduces the approximation on any segment 
$[A,B]$ to the approximation on the segment $[-1,1]$. Therefore, we shall 
consider the case when $A=-1$, $B=1$ in the sequel unless otherwise stated.

For the calculation of integrals, the Gauss--Hermite--Chebyshev quadrature
formula is used:
\begin{equation}
\int\limits_{-1}^1{\varphi (x)\over\sqrt{1-x^2}}\,dx={\pi\over s}
\sum_{i=1}^s\varphi\Big(\cos{2i-1\over 2s}\pi\Big),
\end{equation}
where $s$ is the number of interpolation points;
for polynomials of degree $2s-1$
this formula is exact, so the precision of formula~(19)
increases rapidly as the parameter $s$ increases and depends on the
quality of the approximation of the function $\varphi (s)$ by polynomials.
To calculate the values of Chebyshev polynomials, the well-known 
recurrence relation is applied.

If the function $f(x)$ is even and an approximant is looked for the
form~(17),
then system~(13) is transformed into the following system of equations:
\begin{equation}
\sum_{i=0}^na_i\int\limits_{-1}^1{x^{2i}T_{2k}(x)\over\sqrt{1-x^2}}\,dx
-\sum_{j=
0}^mb_j\int\limits_{-1}^1{x^{2j}T_{2k}(x)f(x)\over\sqrt{1-x^2}}\,dx=0,
\end{equation}
where $k=0,1,\dots ,m+n$. If $f(x)$ is an odd function and an approximant
is looked for in the form~(18), then, first, by means of the solution of
system~(20) complemented by one of the normalization conditions, one
determines an approximant of the form~(17) to the even
function $f(x)/x$, and then the
obtained approximant is multiplied by $x$. This procedure allows us
to avoid a large relative error for $x=0$.

This algorithm is rather simple; for its implementation only two standard
subroutines are needed (for solving systems of linear algebraic equations 
and for numerical integration). However, the algorithm is efficient.

The possibilities of the PADE program are demonstrated in Table~1. This
table contains errors of certain approximants obtained by means of this
program.

\newpage

\centerline {Table 1}
\medskip

\begin{tabular}{|c|c|c|c|c|c|}
\hline
 Function &m&n&$\Delta$&$\delta$&$\delta_{\min}$\\
\hline
$\sqrt{x}$&$2$&$2$&$0.8\cdot 10^{-6}$&$1.13\cdot 10^{-6}$&$0.6\cdot 10^{-6}$\\
\hline
$\sqrt{x}$&$3$&$3$&$1.9\cdot 10^{-9}$&$2.7\cdot 10^{-9}$&$1.12\cdot 10^{-9}$\\
\hline
$\cos{\pi\over 4}x$&$0$&$3$&$0.28\cdot 10^{-7}$&$0.39\cdot 10^{-7}$&$0.32\cdot
10^{-7}$\\
\hline
$\cos{\pi\over 4}x$&$1$&$2$&$0.24\cdot 10^{-7}$&$0.34\cdot 10^{-7}$&$0.29\cdot
10^{-7}$\\
\hline
$\cos{\pi\over 4}x$&$2$&$2$&$0.69\cdot 10^{-10}$&$0.94\cdot 10^{-10}$&$0.79
\cdot 10^{-10}$\\
\hline
$\cos{\pi\over 4}x$&$0$&$5$&$0.57\cdot 10^{-13}$&$0.79\cdot 10^{-13}$&$0.66
\cdot 10^{-13}$\\
\hline
$\cos{\pi\over 4}x$&$2$&$3$&$0.4\cdot 10^{-13}$&$0.55\cdot 10^{-13}$&$0.46
\cdot 10^{-13}$\\
\hline
$\sin{\pi\over 4}x$&$0$&$4$&$0.34\cdot 10^{-11}$&$0.48\cdot 10^{-11}$&$0.47
\cdot 10^{-11}$\\
\hline
$\sin{\pi\over 4}x$&$2$&$2$&$0.32\cdot 10^{-11}$&$0.45\cdot 10^{-11}$&$0.44
\cdot 10^{-11}$\\
\hline
$\sin{\pi\over 4}x$&$0$&$5$&$0.36\cdot 10^{-14}$&$0.55\cdot 10^{-14}$&$0.45
\cdot 10^{-14}$\\
\hline
$\sin{\pi\over 2}x$&$1$&$1$&$0.14\cdot 10^{-3}$&$0.14\cdot 10^{-3}$&$0.12
\cdot 10^{-3}$\\
\hline
$\sin{\pi\over 2}x$&$0$&$4$&$0.67\cdot 10^{-8}$&$0.67\cdot 10^{-8}$&$0.54\cdot
10^{-8}$\\
\hline
$\sin{\pi\over 2}x$&$2$&$2$&$0.63\cdot 10^{-8}$&$0.63\cdot 10^{-8}$&$0.53\cdot
10^{-8}$\\
\hline
$\sin{\pi\over 2}x$&$3$&$3$&$0.63\cdot 10^{-13}$&$0.63\cdot 10^{-13}$&$0.5
\cdot 10^{-13}$\\
\hline
$\tan{\pi\over 4}x$&$1$&$1$&$0.64\cdot 10^{-5}$&$0.64\cdot 10^{-5}$&$0.57\cdot
10^{-5}$\\
\hline
$\tan{\pi\over 4}x$&$2$&$1$&$0.16\cdot 10^{-7}$&$0.16\cdot 10^{-7}$&$0.14\cdot
10^{-7}$\\
\hline
$\tan{\pi\over 4}x$&$2$&$2$&$0.25\cdot 10^{-10}$&$0.25\cdot 10^{-10}$&$0.22
\cdot 10^{-10}$\\
\hline
$\arctan x$&$0$&$7$&$0.75\cdot 10^{-7}$&$10^{-7}$&$10^{-7}$\\
\hline
$\arctan x$&$2$&$3$&$0.16\cdot 10^{-7}$&$0.51\cdot 10^{-7}$&$0.27\cdot 
10^{-7}$\\
\hline
$\arctan x$&$0$&$9$&$0.15\cdot 10^{-8}$&$0.28\cdot 10^{-8}$&$0.23\cdot 
10^{-8}$\\
\hline
$\arctan x$&$3$&$3$&$0.54\cdot 10^{-9}$&$1.9\cdot 10^{-9}$&$0.87\cdot 
10^{-9}$\\
\hline
$\arctan x$&$4$&$4$&$0.12\cdot 10^{-11}$&$0.48\cdot 10^{-11}$&$0.17\cdot 
10^{-11}$\\
\hline
$\arctan x$&$5$&$4$&$0.75\cdot 10^{-13}$&$3.7\cdot 10^{-13}$&$0.71\cdot 
10^{-13}$\\
\hline
\end{tabular}
\medskip

\noindent
For every approximant, the absolute error $\Delta$, the relative
error $\delta$, and (for comparison) the best possible relative
error $\delta_{\min}$ taken from \cite{11} are indicated. The function
$\sqrt{x}$ is approximated on the segment $[1/2,1]$ by the expression of
the form~(1), the function $\cos{\pi\over 4}x$ is approximated on the segment
$[-1,1]$ by the expression of the form~(17), and all the others are
approximated on the same segment by the expression of the form~(18).

\section{Error autocorrection for the PADE program}

The condition numbers of systems of equations that arise while calculating,
by means of the PADE program, the approximants considered above are
very large: for example, while calculating the approximant of the
form~(18) on the segment $[-1,1]$ to $\sin{\pi\over 2}x$ for $m=n=3$, the
corresponding condition number is of order $10^{13}$.
As a result, the coefficients of the approximant are
determined with a large error. In particular, a small perturbation of
the system of linear equations arising when passing from one computer
to another (because of the calculation errors)
gives rise to large perturbations in the coefficients of the approximant.
Fortunately, the effect of error autocorrection improves the situation, 
and the errors of the approximant have no
substantial changes under this perturbation. This fact is
described in the Introduction, where concrete examples are also given.

Consider some more examples connected with passing from one computer
to another (see \cite{6, 7} for details). The branch of the algorithm 
which corresponds
to the normalization condition~(14) (i.e., to $b_0=1$) is considered.
For $\arctan x$ the calculation of an approximant of the form~(18)
on the segment $[-1,1]$ for $m=n=5$ gave an approximant with the 
absolute error $\Delta=0.35\cdot 10^{-12}$
and the relative error $\delta =0.16\cdot 10^{-11}$. The corresponding
system of linear algebraic equations has the condition number of
order $10^{30}$! Passing to another computer we obtain the
following: $\Delta =0.5\cdot 10^{-14}$, $\delta =0.16\cdot 10^{-12}$,
the condition number is of order $10^{14}$, and the errors
$\Delta a_1$ and $\Delta b_1$ in the coefficients $a_1$ and $b_1$ in~(18)
are greater in absolute value than $1$! This example shows that
the problem of computing the condition number of an ill-conditioned system is,
in its turn, ill-conditioned. Indeed, the
condition number is, roughly speaking, determined by values of
coefficients of the inverse matrix, every column of the inverse matrix 
being the solution of the
system of equations with the initial matrix of coefficients, i.e., of
an ill-conditioned system.

Consider in detail the effect of error autocorrection for the
approximant of the form~(17) on the segment $[-1,1]$ to the
function $\cos{\pi\over 4}x$ for $m=2$, $n=3$. 
For our two different computers, two different approximants
were obtained with the coefficients $a_i$,$b_i$ and $\tilde a_i$,$\tilde b_i$
respectively. In the both cases the calculation number is of order $10^9$,
absolute error $\delta = 0.55\cdot10^{-13}$; these errors are closed to 
the best possible. 
The coefficients of the approximants obtained by means of the
computers mentioned above
and the coefficients of the error approximant (see Section 3 above) are
as follows:
$$
\tilde a_0=0.9999999999999600, \qquad a_0=0.9999999999999610, 
$$
$$
\Delta a_0=-10^{-15},
$$
$$
\tilde a_1=-0.2925310453579570, \qquad a_1=-0.2925311264716216, 
$$
$$
\Delta a_1=10^{-7}\cdot 0.811136646,
$$
$$
\tilde a_2=10^{-1}\cdot 0.1105254254716866, \qquad a_2=10^{-1}\cdot 0.1105256585556549, 
$$
$$
\Delta a_2=-10^{-7}\cdot 0.2330839683,
$$
$$
\tilde a_3=10^{-3}\cdot 0.1049474500904401, \qquad a_3=10^{-3}\cdot 0.1049482094850086, 
$$
$$
\Delta a_3=10^{-9}\cdot 0.7593947685,
$$
$$
b_0=1, \qquad \tilde b_0=1, 
$$
$$
\Delta b_0=0,
$$
$$
\tilde b_1=10^{-1}\cdot 0.1589409217324021, \qquad b_1=10^{-1}\cdot 
0.1589401105960337, 
$$
$$
\Delta b_1=10^{-7}\cdot 0.8111363684,
$$
$$
\tilde b_2=10^{-3}\cdot 0.1003359011092697, \qquad b_2=10^{-3}\cdot 0.1003341918083529, 
$$
$$
\Delta b_2=10^{-8}\cdot 0.17093009168.
$$
Thus, the error approximant has the form
\begin{equation}
{\Delta P\over\Delta Q}={\Delta a_0+\Delta a_1x^2+\Delta a_2x^4+\Delta
a_3x^6
\over\Delta b_1x^2+\Delta b_2x^4}.
\end{equation}
If the relatively small quantity $\Delta a_0=-10^{-15}$
in~(21) is omitted,
then, as testing by means of a computer shows (2000 checkpoints), this
expression is an approximant to the function $\cos{\pi\over 4}x$ on the
segment $[-1,1]$ with the absolute and the relative errors
$\Delta =\delta =0.22\cdot 10^{-6}$.

But the polynomial $\Delta Q$ is zero at $x=0$, and the polynomial
$\Delta P$ takes a small, but nonzero value at $x=0$.
Fortunately, relation ~(3) can be rewritten in the following way:
\begin{equation}
{\widetilde{P}\over\widetilde{Q}}-{P\over Q}={\Delta P\over\widetilde{Q}}-
{\Delta Q\over\widetilde{Q}}\cdot{P\over Q}.
\end{equation}
Thus, as $\Delta Q\to 0$, the effect of error autocorrection
arises because the quantity $\Delta P$ is close to
zero, and the error of the approximant $P/Q$ is determined by the
error of the coefficient $a_0$. The same situation also takes place when
the polynomial $\Delta Q$ vanishes at an arbitrary point
$x_0$ belonging to the segment $[A,B]$ where the function is
approximated. It is clear that if one chooses the standard normalization
($b_0=1$), then the error approximant has actually two coefficients less than
the initial one. It is clear that in the
general case the normalization conditions $a_n=1$ or $b_m=1$ result
in the following: the coefficients of the error approximant form an
approximate solution of the homogeneous system of linear algebraic
equations whose exact solution determines the Pad\'e--Chebyshev approximant
having one coefficient less than the initial one. The effect of error
autocorrection improves again the accuracy of this error
approximant; thus,
``the snake bites its own tail''. A situation also
arises in the case when the approximant of the form~(17) to an
even function is constructed by solving the system of equations~(20).

Sometimes it is possible to decrease the error of the approximant by means of
the fortunate choice of the normalization condition. As
an example, consider the approximation of the function $e^x$
on the segment $[-1,1]$ by rational
functions of the form~(1) for $m=15$, $n=0$.
For the traditionally accepted normalization $b_0=1$, the PADE
program yields an
approximant with the absolute error $\Delta =1.4\cdot 10^{-14}$ and
the relative error $\delta =0.53\cdot 10^{-14}$. After passing to the
normalization condition $b_{15}=1$, the errors are reduced nearly one half:
$\Delta =0.73\cdot 10^{-14}$, $\delta =0.27\cdot 10^{-14}$.
Note that the condition number increases: in the first
case it is $2\cdot 10^6$, and in the second case it is $6\cdot 10^{16}$.
Thus the error decreases notwithstanding the fact that the system of
equations becomes drastically ill-conditioned. This example shows that
the increase of accuracy of the error approximant can be
accompanied by the increase
of the condition number, and, as experiments show, by the increase of errors
of the numerator and the denominator of the approximant. The fortunate choice
of the normalization condition depends on the particular situation.

A specific situation arises when the degree of the numerator (or of the
denominator) of the approximant is equal to zero. In this case the unfortunate
choice of the normalization condition results in the following:
the error approximant becomes zero or is not well-defined.
For $n=0$ it is expedient to choose condition~(15), as it was done in the
example given above. For $m=0$ (the case of the polynomial approximation)
it is usually expedient to choose condition~(16). 

One could seek the numerator and the denominator of the approximant in
the form
\begin{equation}
P=\sum_{i=0}^na_iT_i,\quad Q=\sum_{j=0}^mb_jT_j,
\end{equation}
where $T_i$ are the Chebyshev polynomials. In this case the system of linear
equations determining the coefficients would be better
conditioned. But the calculation of the polynomials of the form~(23) by, for
example, the Chenshaw method, results in lengthening the computation time,
although it has a favorable effect upon the error of calculations. The 
transformation of the polynomials
$P$ and $Q$ from the form~(23) into the standard form also 
requires additional efforts.

In practice it is more convenient to use approximants
represented in the form~(1), (17), or~(18), and calculate the fraction's
numerator and denominator according the Horner scheme. In this case
the normalization $a_n=1$ or $b_m=1$ allows to reduce the number of
multiplications. 

The use of the algorithm does not require that
the approximated function be expanded into a series or a continued
fraction beforehand.
Equations~(12) or~(13) and the quadrature formula~(19) show that 
the algorithm uses only the values of the approximated function
$f(x)$ at the interpolation points of the quadrature formula (which are
zeros of some Chebyshev polynomial).

On the segment $[-1,1]$ the linear Pad\'e--Chebyshev approximants give a
considerably smaller error than the classical Pad\'e approximants. For
example, the Pad\'e approximant of the form~(1) to the function $e^x$ for
$m=n=2$ has the absolute error $\Delta (1)=4\cdot 10^{-3}$ at the point
$x=1$, but the PADE program gives an approximant of the same form with
the absolute error $\Delta =1.9\cdot 10^{-4}$ (on the entire the
segment), i.e., the latter is 20 times smaller than the previous one.
The absolute error of the best approximant is $0.87\cdot 10^{-4}$.

\section{The ``cross--multiplied'' linear Pad\'e--Cheb\-yshev approximation}

As a rule, linear Pad\'e--Chebyshev approximants are constructed
according to the following scheme, see, e.g., \cite{21, 11, 16}. Let the
approximated function be decomposed into the series in Chebyshev
polynomials
\begin{equation}
f(x)=\sumpr_{i=0}^{\infty} c_iT_i(x)={1\over
2}c_0+c_1T_1(x)+c_2T_2(x)+\dots,
\end{equation}
where the notation $\sumpr\limits_{i=\infty}^m$ means that the first term $u_0$
in the sum is replaced by $u_0/2$. The rational approximant is
looked for in the form
\begin{equation}
R(x)= \frac{\sumpr\limits_{i=0}^n a_iT_i(x)}
{\sumpr\limits_{j=0}^m b_jT_j(x)}; 
\end{equation}
the coefficients $b_j$ are determined by means of the system of linear
algebraic equations
\begin{equation}
\sumpr_{j=0}^m b_j(c_{i+j}+c_{\vert i-j\vert})=0,\qquad
i=n+1,\dots,n+m, 
\end{equation}
and the coefficients $a_i$ are determined by the equalities
\begin{equation}
a_i={\frac 12}\sumpr_{j=0}^m b_j(c_{i+j}+c_{\vert i-j\vert})=0,\qquad
i=0,1,\dots,n. 
\end{equation}
It is not difficult to verify that this algorithm must lead to the
same results as the algorithm described in Section 5
if the calculation errors are not taken into account.

The coefficients $c_k$ for $k=0,1,\dots,n+2m$, are present in~(26) and~(27),
i.e., it is necessary to have the first $n+2m+1$ terms of series~(24).
The coefficients $c_k$ are known, as a rule, only approximately. To determine
them one can take the truncated expansion of $f(x)$ into the series in
powers of $x$ (the Taylor series) and by means of the well-known economization
procedure transform it into the form
\begin{equation}
\sum_{i=0}^{n+2m} \tilde c_iT_i(x).
\end{equation}

\section{Nonlinear Pad\'e--Chebyshev approximations }

A rational function $R(x)$ of the form~(1) or~(25) is called a
{\it nonlinear Pad\'e--Chebyshev approximant} to the function $f(x)$
on the segment $[-1,1]$, if
\begin{equation}
\int\limits_{-1}^1 \big(f(x)-R(x)\big)T_k(x)w(x)\,dx=0,\qquad
k=0,1,\dots,m+n,
\end{equation}
where $T_k(x)$ are the Chebyshev polynomials, $w(x)=1/\sqrt{1-x^2}$.

In the paper \cite{22} the following algorithm of computing the
coefficients of the approximant indicated above is given. Let the
approximated function $f(x)$ be expanded into series~(24) in
Chebyshev polynomials. Determine the auxiliary quantities $\gamma_i$ from
the system of linear algebraic equations
\begin{equation}
\sum_{j=0}^m\gamma_jc_{\vert k-j\vert}=0,\qquad
k=n+1,n+2,\dots,n+m,
\end{equation}
assuming that $\gamma_0=1$. The coefficients of the denominator in
expression~(25) are determined by the relations
$$
b_j=\mu\sum_{i=0}^{m-j}\gamma_i\gamma_{i+j},
$$
where $\mu^{-1}=1/2\sum_{i=1}^n\gamma_i^2$; this implies
$b_0=2$. Finally, the coefficients of the numerator are determined by
formula~(27). It is possible to solve system~(30) explicitly and to
indicate the
formulas for computing the quantities $\gamma_i$. One can also estimate
explicitly the absolute error of the approximant. This
algorithm is described in detail in the book \cite{20}; see also~\cite{16}.

In contrast to the linear Pad\'e--Chebyshev approximants, the nonlinear
approximants of this type do not always exist, but it is possible to indicate
explicitly verifiable conditions guaranteeing the
existence of such approximants \cite{20}. The nonlinear Pad\'e--Chebyshev
approximants (in comparison with the linear ones) have, as a rule, a somewhat
smaller absolute errors, but can have larger relative errors.
Consider, as an example, the approximant of the form~(1) or~(25)
to the function $e^x$ on the segment $[-1,1]$ for $m=n=3$. In this case
the absolute error for a nonlinear Pad\'e--Chebyshev approximant
is $\Delta =0.258\cdot 10^{-6}$, and the relative error,
$\delta =0.252\cdot 10^{-6}$; for the linear Pad\'e--Chebyshev approximant
$\Delta =0.33\cdot 10^{-6}$ and $\delta =0.20\cdot 10^{-6}$.

\section{Applications of the computer algebra system
REDUCE to the construction of rational approximants}

The computer algebra system REDUCE \cite{23} allows us to handle
formulas at symbolic level and is a convenient tool for the
implementation of algorithms of computing rational approximants.
The use of this system allows us to bypass
the procedure of working out the algorithm of computing the approximated
function if this function is presented in
analytical form or when the Taylor series coefficients are
known or are determined analytically from a differential equation.
The round-off errors can be eli\-minated by using the exact
arithmetic of rational numbers represented in the form of
ratios of integers.

Within the framework of the REDUCE system, the program package for
enhanced precision computations and construction of rational approximants
is implemented; see, for example \cite{4}. In particular, the algorithms
from Sections 6 and 7 (which have similar structure) are
implemented, the approximated function being first expanded into the
power (Taylor) series,
$f=\sum_{k=0}^{\infty}f^{(k)}(0){x^k/k!}$, and then the truncated series
\begin{equation}
\sum_{k=0}^N f^{(k)}(0){x^k\over k!},
\end{equation}
consisting of the first $N+1$ terms of the Taylor series (the value $N$
is determined by the user) being transformed into a polynomial
of the form~(28) by means of the economization procedure.

The algorithms implemented by means of the REDUCE system allow us to obtain
approximants in the form~(1) or~(25), estimates of the absolute and the
relative error, and the error curves. The output includes the Fortran
program of computing the corresponding approximant, the constants of
rational arithmetic being transformed into the standard floating point
form. When computing the values of the obtained approximant, this
approximant can be transformed into the form most convenient for the
user. For example, one can calculate values of the numerator and the
denominator of the fraction of the form~(1) according to the Horner scheme,
and for the fraction of the form~(25), according to Clenshaw's scheme, and
transform the rational expression into a continued fraction or a Jacobi
fraction as well.

The ALGOL-like input language of the REDUCE system and convenient tools for
solving problems of linear algebra guarantee simplicity and
compactness of the programs. For example, the length of the program for
computing linear Pad\'e--Chebyshev approximants is 62 lines.

\section{Error approximants for linear and nonlinear
Pad\'e--Chebyshev approximations }

Relations~(29) can be regarded as a system of equations for
the coefficients of the approximant. Let the approximants
$R(x)=P(x)/Q(x)$ and $\widetilde{R}(x)=\widetilde{P}(x)/\widetilde{Q}(x)$,
where $P(x)$, $\widetilde{P}(x)$ are polynomials of degree $n$ and
$Q(x)$, $\widetilde{Q}(x)$ are polynomials of degree $m$, be obtained
by approximate solving the indicated system of equations. Consider the
error approximant $\Delta P(x)/\Delta Q(x)$, where
$\Delta P(x)=\widetilde{P}(x)-P(x)$, $\Delta Q(x)=\widetilde{Q}(x)-Q(x)$.
Substituting $R(x)$ and $\widetilde{R}(x)$ in~(29) and subtracting one of the
obtained expressions from the other, we see that the following
approximate equality holds:
$$
\int\limits_{-1}^1\bigg({\tilde P(x)\over\tilde Q(x)}-{P(x)\over
Q(x)}\bigg)
T_k(x)w(x)\,dx\approx 0,\qquad k=0,1,\dots,m+n,
$$
which directly implies that $\tilde R(x)-R(x)$ is clause to zero.
This and the equality~(3) imply the approximate equality
\begin{equation}
\int\limits_{-1}^1\bigg({\Delta P(x)\over\Delta Q(x)}-{P(x)\over
Q(x)}\bigg){\Delta Q\over\tilde Q}T_k(x)w(x)\,dx\approx \myone,
\end{equation}
where $k=0,1,\dots,m+n$, $w(x)=1/\sqrt{1-x^2}$. If the quantity
$\Delta Q$ is relatively not small (this is connected with the fact that
the system of equations~(30) is ill-conditioned), then, as follows from
equality~(32), we can naturally expect that the error approximant
is close to $P/Q$ and, consequently, to the approximated function~$f(x)$.

Due to the fact that the arithmetic system of rational numbers is used, the
software described in Section 7 allows us to eliminate the round-off errors 
and to estimate the ``pure'' influence of errors in the approximated function
on the coefficients of the nonlinear
Pad\'e--Chebyshev approximant. In
this case the effect of error autocorrection can be substantiated by a more
accurate reasoning which is valid both for nonlinear Pad\'e--Chebyshev
approximants and for linear ones, and even for the linear generalized
Pad\'e approximants connected with different systems of orthogonal
polynomials. This reasoning is analogous to Y.~L.~Luke's
considerations \cite{13} for the case of classical Pad\'e approximants.

Assume that the function $f(x)$ is expanded into series~(24) and that the
rational approximant $R(x)=P(x)/Q(x)$ is looked for in the form~(25).

Let $\Delta b_j$ be the errors in coefficients of the approximant's
denominator~$Q$. In the linear case these errors arise when solving the
system of equations~(26), and in the nonlinear case, when
solving the system of equations~(27). In both the cases the coefficients
in the approximant's numerator are determined by equations~(27),
whence we have
\begin{equation}
\Delta a_i={\frac 12}{\sumpr_{j=0}^m}\Delta b_j(c_{i+j}+
c_{\vert i-j\vert}),\qquad i=0,1,\dots,n.
\end{equation}
This implies the following fact: the error approximant
$\Delta P/\Delta Q$ satisfies the relations
\begin{equation}
\int\limits_{-1}^1\big(f(x)\Delta Q(x)-\Delta
P(x)\big)T_i(x)w(x)\,dx=0,\qquad i=0,1,\dots,n,
\end{equation}
which are analogous to relations~(12) defining the linear
Pad\'e--Chebyshev approximants. Indeed, let us use the well-known
multiplication formula for Chebyshev polynomials:
\begin{equation}
T_i(x)T_j(x)=\frac 12\big[T_{i+j}(x)+T_{\vert i-j\vert}(x)\big],
\end{equation}
where $i$, $j$ are arbitrary indices; see, for example~\cite{16, 20}. 
Taking~(35) into account, the quantity $f\Delta Q-\Delta P$
can be rewritten in the following way:
$$
f\Delta Q-\Delta P=\bigg({\sumpr_{j=0}^m}\Delta b_jT_j\bigg)\bigg(
{\sumpr_{i=0}^\infty}c_iT_i\bigg)-{\sumpr_{i=0}^n}\Delta a_iT_i
$$
$$
=\frac 12{\sumpr_{i=0}^\infty}\bigg[{\sumpr_{j=0}^m}\Delta b_j(c_{i+j}+
c_{\vert i-j\vert})\bigg]T_i-{\sumpr_{i=0}^n}\Delta a_iT_i.
$$
This formula and~(33) imply that
\begin{equation}
f\Delta Q-\Delta P=O(T_{n+1}),
\end{equation}
i.e., in the expansion of the function $f\Delta Q-\Delta P$ into the
series in Chebyshev polynomials, the first $n+1$ terms {\it are absent},
and the latter is equivalent to relations~(34) by virtue of the fact
that the Chebyshev polynomials form an orthogonal system.

Consider an arbitrary rational function of the form (1) or (8)
$$
R_{m, n}(x)= \frac{a_0+a_1x+\cdots+a_nx^n}{b_0+b_1x+\cdots+b_mx^m}=
\frac{P_n(x)}{Q_m(x)}.
$$
We shall say that $R_{m,n}(x)$ is a {\it generalized linear 
Pad\'e-Chebyshev approximant of order} $N$ to the function $f(x)$ if
$$
\int_{-1}^1\Phi (x)T_k(x)w(x) dx = 0,\quad k=0, 1, \cdots, N,
$$
where $T_k(x)=\cos(n \arccos x)$ are the Chebyshev polynomials,
$w(x)=1/\sqrt{1-x^2}$, $\Phi(x)=f(x)Q_m(x)-P_n(x)$. This means that
the first $N+1$ terms in the expansion of the function
$\Phi(x)$ into the series in Chebyshev polynomials (``the
Fourier-Chebyshev series'') are absent, i.e.
$$
f(x)Q_m(x)-P_n(x)=O(T_{N+1}).
$$
If $N=m+n$, then we have the usual linear Pad\'e-Chebyshev approximant
discussed above in Section 4. Formula (36) means that the following
result is valid.

{\bf Theorem.} {\it
Let $\Delta P\over{\Delta Q}$ be the error approximant to $f(x)$
generated by the approximant~{\rm (25)} for the case of linear or nonlinear
Pad\'e--Chebyshev approximation and algorithms described in Sections {\rm 6}
and {\rm 7}. Then this error approximant $\frac{\Delta P}{\Delta Q}$ is a 
generalized linear Pad\'e-Chebyshev approximant of order $n$ to
the function $f(x)$.}

An equivalent result was discussed in~\cite{5, 7}.
When carrying out actual computations, the coefficients $c_i$ are known
only approximately, and thus the equalities~(33), (34) and (35) are also
satisfied approximately.

\section{Computer experiments for the nonlinear Pad\'e--Chebyshev 
approximation}

Consider the results of computer experiments
that were performed by means of the software implemented
in the framework of the REDUCE system and briefly described
in Section 7 above. At the author's request, computer calculations
were carried out by A.~Ya.~Rodionov. We begin with the example 
considered in Section 5 above,
where the linear Pad\'e--Chebyshev approximant of the form~(17) to the
function $\cos{\pi\over 4}x$ was constructed on the segment
$[-1,1]$ for $m=2$, $n=3$. To construct the corresponding nonlinear
Pad\'e--Chebyshev approximant, it is necessary to specify the value of the
parameter $N$ determining the number of
terms in the truncated Taylor series~(31) of the approximated function.
In this case the  calculation error is determined, in fact, by the
parameter $N$.

The coefficients in approximants of the form~(17) which are obtained for
$N=15$ and $N=20$ (the nonlinear case) and the coefficients in the
error approximant are as follows: 
$$\tilde a_0=0.4960471034987563,\qquad 
a_0=0.4960471027757504,$$
$$\Delta a_0=10^{-8}\cdot 0.07230059,$$
$$ \tilde a_1=-0.1451091945278387,\qquad 
a_1=-0.1451091928755736,$$
$$\Delta a_1=-10^{-8}\cdot 0.16522651,$$
$$\tilde a_2=10^{-2}\cdot 0.5482586543334515,\qquad
a_2=10^{-2}\cdot 0.548258121085953,$$
$$\Delta a_2=-10^{-9}\cdot 0.42224856,$$
$$\tilde a_3=-10^{-4}\cdot 0.5205903601778259,\qquad 
a_3=-10^{-4}\cdot 0.5205902238186334,$$
$$\Delta a_3=-10^{-10}\cdot 0.13635919,$$
$$\tilde b_0=0.4960471034987759,\qquad 
b_0=0.4960471027757698, $$
$$\Delta b_0=10^{-8}\cdot 0.07230061,$$
$$\tilde b_1=10^{-2}\cdot 0.7884201590727615,\qquad 
b_1=10^{-2}\cdot 0.7884203019999351, $$
$$\Delta b_1=-10^{-10}\cdot 0.1429272,$$
$$\tilde b_2=10^{-4}\cdot 0.4977097973870693,\qquad 
b_2=10^{-4}\cdot 0.4977100977750249,$$
$$\Delta b_2=-10^{-10}\cdot 0.300388.$$

Both the approximants have absolute errors $\Delta$ equal to
$0.4\cdot 10^{-13}$
and the relative errors $\delta$ equal to $0.6\cdot 10^{-13}$, these
values being close to the best possible.
The condition number of the system of equations~(30) in
both the cases is $0.4\cdot 10^8$. The denominator $\Delta Q$ of the error
approximant is zero for $x=x_0\approx 0.70752\dots$; the
point $x_0$ is also close to the root of the numerator $\Delta P$ which for
$x=x_0$ is of order $10^{-8}$. Such a situation was
considered in Section 5 above. Outside a small neighborhood of
the point $x_0$ the absolute and the relative errors have the
same order as in the ``linear case'' considered in Section 5. 

Now consider the nonlinear Pad\'e--Chebyshev approximant of the form (17)
on the segment $[-1,1]$ to the function $\tan{\pi\over 4}x$ for
$m=n=3$. In this case the Taylor series converges very slowly, and,
as the parameter $N$ increases, the values of coefficients of the rational
approximant
undergo substantial (even in the first decimal digits) and intricate changes.
The situation is illustrated in Table~2,
where the following values are given: the absolute errors $\Delta$,
the absolute errors $\Delta_0$ of error approximants
(there the approximants are compared for $N=15$ and $N=20$,
for $N=25$ and $N=35$,
for $N=40$ and $N=50$), and also the values of the condition number
$cond$ of the system of linear algebraic equations~(30). In this case
the relative errors coincide with the absolute ones. The best possible
error is $\Delta_{\min}=0.83\cdot 10^{-17}$. 
A small neighborhood of the root of the polynomial $\Delta Q$ is 
eliminated as before. 

\centerline{Table 2}
\medskip

{\footnotesize
\noindent
\begin{tabular}{|c|c|c|c|c|c|c|c|}
\hline
N&$15$&$20$&$25$&$35$&$40$&$50$\\
\hline 
$cond$&$0.76\cdot 10^7$&$0.95\cdot 10^8$&$0.36\cdot 10^{10}$&
$0.12\cdot 10^{12}$&$0.11\cdot 10^{12}$&$0.11\cdot 10^{12}$\\
\hline 
$\Delta$&$0.13\cdot 10^{-4}$&$0.81\cdot 10^{-6}$&$0.13\cdot 10^{-7}$&
$0.12\cdot 10^{-10}$&$0.75\cdot 10^{-12}$&$0.73\cdot 10^{-15}$\\
\hline
$\Delta_0 $&
$ 0.7\cdot 10^{-4}$&&
$0.7\cdot 10^{-8}$&&
$0.2\cdot 10^{-9}$&\\
\hline 
\end{tabular}
}

\section{Small deformations of approximated functions and
acceleration of convergence of series}

Let a function $f(x)$ be expanded into the series in Chebyshev
polynomials, i.e., suppose that $f(x)=\sum_{i=0}^\infty c_iT_i$; 
consider a partial sum
\begin{equation}
\hat f_N(x)=\sum_{i=o}^Nc_iT_i
\end{equation}
of this series. Using formula~(35), it is easy to verify that the
linear Pad\'e--Chebyshev
approximant of the form~(1) or~(25) to the function $f(x)$ coincides
with
the linear Pad\'e--Chebyshev approximant to polynomial~(37) for
$N=n+2m$, i.e., it depends only on the first $n+2m+1$ terms of the
Fourier--Chebyshev series of the function $f(x)$; a similar result is
valid for the approximant of the form~(17) or~(18) to even or odd
functions, respectively. Note that for $N=n+2m$ the polynomial
$\hat f_N$
is the result of application of the algorithm of linear (or nonlinear)
Pad\'e--Chebyshev approximation to $f(x)$, the exponents $m$ and $n$
being replaced by $0$ and $2m+n$.

The interesting effect mentioned in \cite{2} consists in the fact that
the error
of the polynomial approximant $\hat f_{n+2m}$ depending on $n+2m+1$
parameters can exceed the error of the corresponding Pad\'e--Chebyshev
approximant of the form~(1) which depends only on $n+m+1$ parameters.
For example, consider an approximant of the form~(18) to the function
$\tan{\pi\over 4}x$ on the segment $[-1,1]$. For $m=n=3$ the linear
Pad\'e--Chebyshev approximant to $\tan{\pi\over 4}x$ has the error of
order $10^{-17}$, and the corresponding polynomial approximant of the
form~(37) has the error of order $10^{-11}$. This polynomial of degree 19
(odd functions are in question, and hence
$m=n=3$ in~(18) corresponds to $m=6$, $n=7$ in~(1))
can be regarded as a result of a deformation of the approximated function
$\tan{\pi\over 4}x$. This deformation does not affect the first
twenty terms in the expansion of this function in Chebyshev polynomials
and, consequently, does not affect the coefficients in the corresponding
rational Pad\'e--Chebyshev approximant, but leads to an
increase of several orders in its error. Thus, {\it a small deformation of 
the approximated function can
result in a sharp change in the order of error of a rational approximant}.

Moreover the effect just mentioned means that the algorithm extracts an 
additional information concerning the next components of the
Fourier--Chebyshev series from the polynomial~(37). In other words, in 
this case the transition from the Fourier--Chebyshev 
series to the corresponding Pad\'e--Chebyshev approximant accelerates
convergence of the series. A similar effect of acceleration of
convergence of power series by passing to the classical Pad\'e approximant is
known (see, e.g.,~\cite{16}).

It is easy to see that the nonlinear Pad\'e--Chebyshev approximant
of the form~(1) to the function $f(x)$ depends only on the first
$m+n+1$ terms of the Fourier--Chebyshev series for $f(x)$, so that for such
approximants a more pronounced effect of the type indicated above
takes place.

Since one can change the ``tail'' of the Fourier--Chebyshev series in
a quite arbitrary way without affecting the rational Pad\'e--Chebyshev
approximant, the effect of acceleration of convergence can take
place only for the series with an especially regular behavior (and
for the corresponding ``well-behaved'' functions). 
See \cite{5, 7} for some details. 

\section{Error autocorrection and Interval \newline Analysis}

Undoubtedly one of the most actual problems in Interval Analysis
in the sense of \cite{8}--\cite{10} is to get realistic interval estimates 
for calculation errors, i.e. to get efficient estimates close to the virtual 
calculation errors. Difficulties arise where intermediate errors compensate
each other.
For the sake of simplicity let us suppose that we calculate
values of a real smooth function $z=F(y_1,\dots ,y_N)$ of 
real variables $ y_1,\dots ,y_N $, and suppose that all round-off
errors are negligible with respect to input data errors. This
situation is especially examined in the framework of Ju.V. Matijasevich's
``a posteriori interval analysis'', see, e.g., \cite{10}. 
In this case the error
$\Delta F $ of $ F(y_1,\dots ,y_N)$ can be evaluated by the formula
\begin{equation}
\Delta F=\sum^N_{i=1} \frac{\partial F}{\partial y_i}\cdot \Delta y_i+r,
\end{equation}
where r is negligible. The sum 
$\sum^N_{i=1}(\partial F/\partial{y_i})\cdot \Delta y_i$ 
in~(38) can be treated as a scalar product of the gradient
vector $\{\partial F/\partial{y_i}\}$ and the vector of errors 
$\{\Delta y_i\}$.

The effect of error autocorrection corresponds to the case, where 
the gradient $\{\partial F/\partial{y_i}\}$ is not small but the 
scalar product
is small enough. In this case these vectors are almost orthogonal
and the following approximate equation holds:
\begin{equation}
\sum^N_{i=1}\pd{F}{y_i}\cdot \Delta y_i\approx 0
\end{equation}
This effect is typical for some ill-posed problems.

For all the standard interval methods, the best possible estimation 
for $\Delta F$ is given by the formula
\begin{equation}
\vert \Delta F \vert \le \sum^N_{i=1}\vert\frac {\partial F}{\partial y_i}
\vert\cdot \vert \Delta y_i \vert. 
\end{equation}
This estimate is good if the errors $\Delta y_i$ are ``independent''
but it is not realistic in the case discussed
in this paper (calculation of values of rational approximants when the error
autocorrection effect works). In this case 
$$ F(y_1,\dots ,y_N)= R(x,a_0,
\dots,a_n,b_0,\dots,b_m)={a_0+a_1x+a_2x^2+\dots +a_nx^n\over{b_0+b_1x+b_2x^2
+\dots +b_mx^m}},
$$
 where $N=m+n+3$, $ \{y_1,\dots,y_N\}=\{ x,a_0, \dots,a_n,b_0,\dots,b_m \}$.
 For the sake of simplicity let us suppose that 
$\Delta x=0$. In this case we can use the equality~(22) to 
transform the formula~(38) into the formula
\begin{equation}
\Delta R \approx \frac{\Delta P}{Q}-\frac{P\Delta Q}{Q^2}=
\sum^n_{i=0}\frac{x^i}{Q(x)}\Delta a_i+\sum^m_{j=0}
\frac{P(x)x^j}{Q^2(x)}\Delta b_j. 
\end{equation}
So the estimation (39) transforms into the estimation 
\begin{equation}
\Delta R \le \sum^n_{i=0}{\vert x^i \vert \over{\vert Q(x)\vert}}\cdot 
\vert \Delta a_i\vert +\sum^m_{j=0}{\vert P(x)x^j\vert\over{Q^2(x)}}\cdot 
\vert \Delta b_j \vert. 
\end{equation}
It is easy to check that estimations of this type are not
realistic. Consider the following example discussed in Introduction:
$f(x)=\arctan x$ on the segment $[-1,1]$, $R(x)$ has the form (1)
for $m=8$, $n=9$. In this case the estimation (41) is of order 
$10^{-4}$ but in fact $\Delta R$ is of order $10^{-11}$.
This situation is typical for examples examined in this paper.
In fact we have an approximate equation
\begin{equation}
\Delta R\approx\sum^n_{i=0} \frac{x^i}{Q(x)}\Delta a_i+\sum^m_{j=0}
\frac{P(x)x^j}{Q^2(x)}\Delta b_j\approx\varepsilon\approx 0,
\end{equation}
where $\varepsilon$ is the absolute error of the used approximation
method. Of course, this approximate equation holds if our
approximation method is good and the uncertainty relation
(discussed above in Section 2) is valid. Then approximate 
equation (43) corresponds to the approximate equation (39).

The error autocorrection effect appears not only in rational approximation
but in many problems. Other examples (where this effect occurs in a 
weaker form) are the method of least squares and some popular methods
for the numerical integration of ordinary and partial differential 
equations, see, e.g., \cite{24}--\cite{29}.   

In principle the error autocorrection effect appears if input
data form an (approximate) solution (or solutions) of an equation
(or equations or systems of equations). Then the corresponding
errors could form an approximate solution (or solutions) for
another equation (or equations or systems of equations). As a result
this could lead to corrections for standard interval estimates.

Of course, in the theory we can include all the preliminary
numerical problems to our concrete problem and to use, e.g., a
posteriori interval analysis for the ``united'' problem.
However, in practice this is not convenient.

In practice situations of this kind often appear if we use
approximate solutions to ill conditioned systems of linear 
algebraic equations. If the condition number of the system is great
and the residual of the solution (with respect to the system)
is small, then our software must send us a ``warning''.
This means that an additional investigation for error estimates is
needed. In the theory of interval analysis this corresponds
to a further development of ``a posteriori interval methods'' in the
spirit of \cite{10}, \cite{27}--\cite{30}.

{\bf Remark.} We have  discussed ``smooth'' computations. Note
that for many ``nonsmooth'' optimization problems all the interval
estimates could be good and absolutely exact. A situation of
this kind (related to solving systems of linear algebraic
equations over idempotent semirings) is described in \cite{31}.
 

%
\end{document}